\documentclass[12pt,leqno]{amsart}

\usepackage{amsmath, amssymb,latexsym}
\usepackage{amsthm,amscd,amsfonts}
\usepackage{graphicx} 

\begin{document}




\newtheorem{thm}{Theorem}[section]
\newtheorem{prop}[thm]{Proposition}
\newtheorem{lem}[thm]{Lemma}
\newtheorem{cor}[thm]{Corollary}
\newtheorem{conj}[thm]{Conjecture}
\newtheorem{ddef}[thm]{Definition}
\newtheorem{ex}[thm]{Example}
\newtheorem{rem}[thm]{Remark}
\newtheorem{notation}[thm]{Notation}

\numberwithin{equation}{section}
\numberwithin{figure}{section}

\newcommand{\bthm}{\begin{thm}} \newcommand{\ethm}{\end{thm}}
\newcommand{\bthms}{\begin{thm*}} \newcommand{\ethms}{\end{thm*}}
\newcommand{\blem}{\begin{lem}} \newcommand{\elem}{\end{lem}}
\newcommand{\bcor}{\begin{cor}} \newcommand{\ecor}{\end{cor}}
\newcommand{\bprop}{\begin{prop}} \newcommand{\eprop}{\end{prop}}
\newcommand{\bproof}{\begin{proof}} \newcommand{\eproof}{\end{proof}}
\newcommand{\bddef}{\begin{ddef}} \newcommand{\eddef}{\end{ddef}}
\newcommand{\bconj}{\begin{conj}} \newcommand{\econj}{\end{conj}}
\newcommand{\brem}{\begin{rem}} \newcommand{\erem}{\end{rem}}
\newcommand{\bca}{\begin{cases}} \newcommand{\eca}{\end{cases}}

\newcommand{\att}[1]{\tcr{{\em #1}}}

\newcommand{\beq}{\begin{equation}} \newcommand{\eeq}{\end{equation}}
\newcommand{\beqs}{\begin{equation*}} \newcommand{\eeqs}{\end{equation*}}
\newcommand{\beqa}{\begin{eqnarray}} \newcommand{\eeqa}{\end{eqnarray}}
\newcommand{\beqas}{\begin{eqnarray*}} \newcommand{\eeqas}{\end{eqnarray*}}
\newcommand{\barr}{\begin{array}} \newcommand{\earr}{\end{array}}
\newcommand{\btab}{\begin{tabular}} \newcommand{\etab}{\end{tabular}}

\newcommand{\bit}{\begin{itemize}} \newcommand{\eit}{\end{itemize}}
\newcommand{\ben}{\begin{enumerate}} \newcommand{\een}{\end{enumerate}}
\newcommand{\bce}{\begin{center}} \newcommand{\ece}{\end{center}}

\newcommand{\defeq}{\stackrel{\rm def}{=}}
\newcommand{\noset}{\varnothing}
\newcommand{\bd}{\partial}

\newcommand{\ho}{\hat{1}}
\newcommand{\hz}{\hat{0}}
\newcommand{\whz}{\widehat{0}}
\newcommand{\who}{\widehat{1}}
\newcommand{\wh}{\widehat}
\newcommand{\wt}{\widetilde}
\newcommand{\cover}{\prec}
\newcommand{\cov}{\lhd}
\newcommand{\cocov}{\rhd}
\newcommand{\opp}{\mathrm{opp}}
\newcommand{\join}{\vee}
\newcommand{\meet}{\wedge}

 
\renewcommand{\AA}{\mathcal{A}}
\newcommand{\BB}{\mathcal{B}}
\newcommand{\CC}{\mathcal{C}}
\newcommand{\DD}{\mathcal{D}}
\newcommand{\EE}{\mathcal{E}}
\newcommand{\FF}{\mathcal{F}}
\newcommand{\HH}{\mathcal{H}}
\newcommand{\II}{\mathcal{I}}
\newcommand{\LL}{\mathcal{L}}
\newcommand{\MM}{\mathcal{M}}
\newcommand{\NN}{\mathcal{N}}
\newcommand{\PP}{{\mathcal{P}}}
\newcommand{\NP}{\mathcal{NP}}
\newcommand{\QQ}{\mathcal{Q}}
\newcommand{\RR}{\mathcal{R}}
\renewcommand{\SS}{\mathcal{S}}


\newcommand{\C}{\mathbb{C}}
\newcommand{\F}{\mathbb{F}}
\renewcommand{\L}{\mathbb{L}}
\newcommand{\N}{\mathbb{N}}
\renewcommand{\P}{\mathbb{P}}
\newcommand{\Q}{\mathbb{Q}}
\newcommand{\R}{\mathbb{R}}
\newcommand{\Z}{\mathbb{Z}}


\newcommand{\al}{\alpha}
\newcommand{\be}{\beta}
\newcommand{\de}{\delta} \newcommand{\De}{\Delta} 
\newcommand{\ga}{\gamma} \newcommand{\Ga}{\Gamma}
\newcommand{\la}{\lambda} \newcommand{\La}{\Lambda}
\newcommand{\om}{\omega} \newcommand{\Om}{\Omega}
\newcommand{\si}{\sigma} \newcommand{\Si}{\Sigma}
\newcommand{\ze}{\zeta}
\newcommand{\vphi}{\varphi}
\newcommand{\eps}{\varepsilon}

\newcommand{\X}{\bf{X}}
\newcommand{\norm}[1]{\lVert#1\rVert}
\newcommand{\abs}[1]{\lvert#1\rvert}
\newcommand{\qbinom}[2]{\left[\ba{c}{#1}\\{#2}\ea\right]} 
\newcommand{\comp}{\models}

\newcommand{\vanish}[1]{}

\newcommand{\st}{\,:\,} 
\newcommand{\sbseq}{\subseteq}
\newcommand{\spseq}{\supseteq}


\newcommand{\larr}{\leftarrow}
\newcommand{\rarr}{\rightarrow}
\newcommand{\Larr}{\Leftarrow}
\newcommand{\Rarr}{\Rightarrow}
\newcommand{\lrarr}{\leftrightarrow}
\newcommand{\Lrarr}{\Leftrightarrow}

\newcommand{\longlarr}{\longleftarrow}
\newcommand{\longrarr}{\longrightarrow}
\newcommand{\Longlarr}{\Longleftarrow}
\newcommand{\Longrarr}{\Longrightarrow}
\newcommand{\longlrarr}{\longleftrightarrow}
\newcommand{\Longlrarr}{\Longleftrightarrow}

\title[A cell complex in number theory]
{A cell complex in number theory}
\author[Anders Bj\"orner]{Anders Bj\"orner   \\ \\ \\
{\em DEDICATED TO DENNIS STANTON ON THE OCCASION OF HIS 60TH BIRTHDAY}}

\address{Institut Mittag-Leffler, Aurav\"agen 17, S-182 60 Djursholm, Sweden}
\email{bjorner@mittag-leffler.se}
\address{Kungl. Tekniska H\"ogskolan, Matematiska Inst.,
 S-100 44 Stockholm, Sweden}
\email{bjorner@math.kth.se}
\thanks{Research supported by the Knut and Alice Wallenberg Foundation,
grant KAW.2005.0098}

\newcommand{\tDn}{\wt{\De}_n}
\newcommand{\sqf}{squarefree}

\begin{abstract} 
Let $\De_n$ be the simplicial complex of squarefree positive integers
less than or equal to $n$ ordered by divisibility.
It is known that
the asymptotic rate of growth of its Euler characteristic
(the Mertens function) is closely related to deep properties of
the prime number system.

In this paper we study the asymptotic behavior of the individual Betti numbers
$\be_k(\De_n)$ and of their sum. 
We show that $\De_n$ has the homotopy type of a wedge of spheres, and that
as $n\rarr\infty$
$$\sum \be_k(\De_n)  = \frac{2n}{\pi^2} + O(n^{\theta}),\;\;  \mbox{ for all } \theta > \frac{17}{54}.$$
Furthermore, for fixed $k$, 
$$\be_k(\De_n) \sim \frac{n}{2 \log n} \frac{(\log\log n)^k}{k!}$$
 As a number-theoretic byproduct we obtain inequalities
$$\partial_{k}\left(\si_{k+1}^{odd}(n)\right) \le \si_{k}^{odd}(n/2),$$
where $\si_k^{odd}(n) $ denotes the number of odd \sqf\ integers $\le n$ with
 $k$ prime factors, and $\partial_{k}$ is a certain combinatorial shadow function.

We also study a CW complex  $\tDn$ that extends the previous simplicial complex.
In $\tDn$ {\em all} numbers $\le n$ correspond to cells
and its Euler characteristic is the summatory Liouville function.
This cell complex $\tDn$ is shown to be homotopy equivalent to a wedge of spheres, and
as $n\rarr\infty$
$$\sum \be_k(\tDn)  = \frac{n}{3} + O(n^{\theta}),\;\;  \mbox{ for all } \theta > \frac{22}{27}.$$

\end{abstract}

\maketitle

\section{Introduction}

Let $M(n)\defeq\sum_{k=1}^n \mu(k)$, where $\mu(k)$ is the number-theoretic M\"obius function.
The rate of growth of the function $M(n)$ is of great interest and  importance in number theory,
as is clear from the following facts:
$$ 
\begin{array}{lll}
\mbox{Prime Number Theorem} & \Longleftrightarrow & 
|M(n)| \leq \varepsilon n,  
  \mbox{    for all } \varepsilon>0 \\
&& \mbox{ and all sufficiently large $n$. } \\
\mbox{Riemann Hypothesis} & \Longleftrightarrow &
|M(n)| \leq n^{1/2+\varepsilon},  
  \mbox{    for all } \varepsilon>0 \\
&& \mbox{ and all sufficiently large $n$. }
\end{array}
  $$
The key  questions concerning  the growth of $M(n)$ are subtle.
For instance, Mertens conjectured in 1897 that
$|M(n)| \leq n^{1/2}  \mbox{    for all sufficiently large $n$. }$
This conjecture was disproved in 1985 by Odlyzko and te Riele.
See e.g. the books by Hardy and Wright \cite{HW} and Ivi\'c  \cite[especially  \S 1.9]{Iv}
for more information about these matters.
  \smallskip
  
This paper has its  genesis in the observation that the Mertens function $M(n)$ can be
interpreted as the Euler characteristic of a simplicial complex. Namely, for each
positive \sqf\  integer $k$, 
let $P(k)$ be the set of its prime factors. 
For instance, $P(165) =\{3,5,11\}$. Then
the set family $$\De_n\defeq \{ P(k)\ :\  \mbox{$k$ is \sqf\ and $k\le n$}\}$$
is closed under taking subsets. In other words, it is an
abstract  simplicial complex. (Remark:
The integer $1$ is \sqf, so we include $P(1)=\varnothing$ in $\De_n$.)
%



 Since, by definition,
 $$ \mu(k)=
 \bca
(-1)^{|P(k)|}, \mbox{ if $k$ is \sqf\ } \\
0, \mbox{ otherwise}
 \eca
 $$
 it follows that 
\beq M(n) = - \chi(\De_n), \label{eulerchar1}
\eeq
where $\chi(\De_n)$ is the reduced Euler characteristic 
(the ordinary Euler characteristic minus one).
Let $\be_k(\De_n)$ denote the $k$-th Betti numbers of reduced simplicial homology, 
i.e., $\be_k(\De_n)\defeq \mbox{rank } \widetilde{H}_k(\De_n,\Z)$. Then, from
the Euler-Poincar\'e formula and equation (\ref{eulerchar1}) we have
\beq M(n) = \sum_{k\ge 0} (-1)^{k-1}\be_k(\De_n). \label{eulerchar2}
\eeq

Thus, important number-theoretic propositions,
such as the Prime Number Theorem and the Riemann Hypothesis, are equivalent to
statements about the asymptotic rate of growth of the Euler characteristic $\chi(\De_n)$.
It therefore seems reasonable to inquire about the Betti numbers of the complex $\De_n$
and their asymptotics.

We prove the following estimates: As $n\rarr\infty$

\beqas
\sum_{k\ge 0} \be_k (\De_n) &=& \frac{2n}{\pi^2} +O(n^{\theta}),\;\;  \mbox{ for all } \theta 
> \frac{17}{54}
 \label{sumall}\\
\sum_{k\ \mathrm{even}} \be_k (\De_n) &\sim& \frac{n}{\pi^2}  \\ 
\sum_{k\ \mathrm{odd}} \be_k (\De_n)& \sim& \frac{n}{\pi^2} \\
\mbox{For fixed $k$: }\; \be_k(\De_n) &\sim& \frac{n}{2 \log n} \frac{(\log\log n)^k}{k!} 
\eeqas
Unfortunately these results fall short of shedding new light on the rate of growth of
the Euler characteristic $M(n)$. Perhaps a study of deeper topological
invariants of $\De_n$ could add something of value.

Let $\si_k^{odd}(n) $ be the number of odd \sqf\ integers $\le n$ with
 $k$ prime factors. 
It turns out that the homology information of $\De_n$ can be expressed in terms of 
these numbers, which, in turn, leads to some purely number-theoretic consequences.
Namely, for certain ``shadow functions'' $\partial_k$ and $\partial^k$,
well-known in extremal combinatorics and
defined in Section 4, we obtain inequalities
\ben
\vspace{3mm}
\item  $\partial_{k}\left(\si_{k+1}^{odd}(n)\right) \le \si_{k}^{odd}(n/2)   $
\vspace{5mm}
\item  $\partial^{k}\left(\si_{2k+2}^{odd}(n)+\si_{2k+1}^{odd}(n)\right) \le \si_{2k}^{odd}(n/2)+\si_{2k-1}^{odd}(n/2)   $
\een

In Section 5 we study a CW complex  $\tDn$ that extends the previous simplicial complex.
In $\tDn$ {\em all} numbers $\le n$ correspond to cells and its Euler characteristic is
the summatory Liouville function. We get a sequence of embeddings
 $$\wt{\De}_2 \hookrightarrow \wt{\De}_3  \hookrightarrow  \cdots \hookrightarrow 
 \wt{\De}_n \hookrightarrow \cdots
$$
that can be seen as a filtration of  the join of denumerably many copies of certain infinite-dimensional
spaces such as $S^{\infty}$ or $\R P^{\infty}$.
The cell complex $\tDn$ is shown to be homotopy equivalent to a wedge of spheres, and
as $n\rarr\infty$
$$\sum \be_k(\tDn)  = \frac{n}{3} + O(n^{\theta}),\;\;  \mbox{ for all } \theta > \frac{22}{27}.$$

The construction of $\tDn$ is based on a more general construction for multicomplexes
given in \cite{BV}. In the last section we recall and expand on
 the details of the construction from \cite{BV}.
 
 Helpful remarks from anonymous referees are gratefully acknowledged.
 
\section{Preliminaries}

Later on we need to refer to a few results from combinatorial topology and number theory.
Here this information is recalled.

A simplicial complex $\De$ on a linearly ordered vertex set 
$x_1, x_2, \ldots, x_t$ is said to be {\em shifted} if for every $F\in \De$:
$$\mbox{ $i<j$, $x_i\notin F$ and $x_j\in F \;\Rightarrow\;  F\setminus \{x_j\} \cup \{x_i\}\in \De$}$$

\bthm \cite[p. 292]{BK} \label{shift}
Suppose that $\De$ is a shifted complex on the vertex set $x_1, x_2. \ldots, x_t$.
Then $\De$ has the homotopy type of a wedge of spheres. Its Betti numbers are
$$\be_k(\De)= \mbox{number of $F\in\De$ such that $ F\cup \{x_1\}\notin \De$ and $|F|=k+1.$}
$$
\ethm

The number $\Omega(n)$ of prime factors
of a positive integer $n$ is called its {\em weight}. I.e.,  by definition, 
$$\Om (p_{i_1}^{e_1}p_{i_2}^{e_2}\cdots  p_{i_k}^{e_k})=e_1 + e_2 + \cdots + e_k.
$$
It is convenient to define the following integer-valued functions for all
positive real numbers  $x$ (not only for integers).
\bddef
For $x\in \R^+$, let
\begin{enumerate}
\item $\si(x) \defeq $ the number of \sqf\ integers in $(0,x]$, 
\item $\si^{odd}(x) \defeq $ the number of odd \sqf\ integers in $(0,x]$, 
\item $\si_k(x) \defeq $ the number of \sqf\ integers in $(0,x]$ of weight $k$,
\item $\si_k^{odd}(x) \defeq $ the number of odd \sqf\ integers in $(0,x]$ of weight $k$,
\item and analogously $\si^{even}(x)$ and  $\si_k^{even}(x)$.
\end{enumerate}
\eddef
\noindent
Note that  $\si(x)=0$ for $x<1$ and similarly for the other functions.
\medskip

The following estimates belong to the classics of number theory.
\bthm \label{classics} As $x\rarr \infty$ we have that
\bit
\item[(a)]  (Gegenbauer, 1885; Jia, 1993; see \cite[pp. 355,359]{HW})
$$ \si(x) = \frac{6x}{\pi^2} +O(x^{\theta}),\;\;  \mbox{ for all } \theta > \frac{17}{54}$$
\item[(b)]  (Landau, 1900; see \cite[p. 491]{HW}) For fixed $k$,
$$\si_k(x) \sim \frac{x}{ \log x} \frac{(\log\log x)^{k-1}}{(k-1)!}$$
\eit
\ethm

\brem {\rm(i) 
Gegenbauer's estimate of the error term was $O(\sqrt{n})$. The sharper
exponent cited here is due to Jia.

(ii)  Landau's asymptotic formula for $\si_k(x)$
was conjectured by Gauss. 
Note that the $k=1$
case is the Prime Number Theorem. Estimates of the error term exist but will not be used
in this paper.
}\erem

\section{The number-theoretic simplicial complex $\De_n$}

Some of the basic concepts of elementary number theory have direct geometric meaning
for the complex $\De_n$. In fact, one can read parts of the book by Hardy and Wright \cite{HW}
as describing the size and shape of this complex.

 For instance,  
the dimension of a simplex $P(k)$ is the weight of $k$ minus one:
$$ \dim P(k) = \Omega(k)-1 $$
and the typical dimension of a simplex in $\De_n$ is circa $\lfloor \log\log n \rfloor$.
For integers $k\ge 1$, let $k\#!$ denote the {\em primorial number}
$p_1 p_2 \cdots p_k$, that is, the product of the first $k$ prime numbers.
Then
$$\dim \De_n= \max \{k\ :\  k\#! \le n\}-1$$
The {\em round numbers} \cite[p. 476]{HW} less than $n$
  (roughly) correspond to the high-dimensional simplices in $\De_n$.
 
 For example, $\dim (\De_{10^7})=7$ and a 
 typical simplex of $\De_{10^7}$ is $2$-dimensional,
  whereas $\dim (\De_{10^{80}})=44$  
 and a typical simplex is of dimension  $4$ or $5$.
Compare this to the discussion in  \cite[p. 477]{HW}, where in passing it is mentioned
that $10^{80}$ is roughly the number of protons in the universe.

\bthm\label{betti1} The complex $\De_n$ is shifted. Its Betti numbers are
\beq
\be_k (\De_n) = \si_{k+1}^{odd} (n) - \si_{k+1}^{odd} (n/2) \eeq\label{betti.shift}
\ethm
\bproof
The vertices of $\De_n$ are the  prime numbers $p_1=2, p_2=3, p_3=5, \ldots$
If for a \sqf\ number $k\le n$ one of its prime factors is replaced by a smaller prime
number, then the new number $k'$ so obtained satisfies $k'<k$, and in particular $k'<n$.
Hence $\De_n$ is shifted.

According to Theorem \ref{shift} the Betti numbers are
$$\be_k(\De_n)= \#\{\mbox{odd \sqf\ integers $b\le n$ such that  $2b \not\le n$
 and $\Omega(b)=k+1$} \},$$
which agrees with the stated expression.
\eproof

The plan from here on is to get rid of the "odd" condition in
formula (\ref{betti.shift}). 
We seek to instead
express the Betti numbers in terms of the standard number theoretic
functions $\si(x)$ and  $\si_k(x)$, so that we can
benefit from the known estimates for these functions.

\blem\label{lem1}
$$ \si_{k}^{even}(x)= \si_{k-1}^{odd}(x/2) 
$$
\elem
\bproof
Multiplication by $2$ gives a bijective map \\[3mm]
\hspace*{4mm} $\{\mbox{odd \sqf\ numbers  in $(0,x/2]$ of weight $k-1$}\} \\
 \hspace*{34mm} \leftrightarrow \{\mbox{even \sqf\ numbers in $(0,x]$ of weight $k$}\}.$
\eproof

\blem \label{lem2}
\beqa
 \si^{odd}_k(x)  &=& \si_k(x)- \si_{k-1}(x/2)+\si_{k-2}(x/4)-\si_{k-3}(x/8)+ \cdots  \label{eq1} \\
\si^{odd}(x)  &=& \si(x)- \si(x/2)+\si(x/4)-\si(x/8)+ \cdots  \label{eq2}
\eeqa
\elem
\bproof
Lemma \ref{lem1} implies that 
$$\si_k (x) = \si^{odd}_k (x) + \si^{odd}_{k-1}(x/2) ,
$$
from which equation (\ref{eq1}) follows. Equation (\ref{eq2})
is then obtained  by summation over all $k$.
\eproof


\bthm\label{bettithm1} As $n\rarr \infty$ we have that
$$ \sum_{k\ge 0} \be_k (\De_n) = \frac{2n}{\pi^2} +O(n^{\theta}),\;\;  \mbox{ for all } \theta > \frac{17}{54}$$
\ethm

\bproof 
Let $N\defeq \lceil \log_2(n)  \rceil$.
We then have that $n\le 2^N < 2n$. Using Theorem \ref{classics}
and Lemma  \ref{lem2} we compute
\beqas
\si^{odd}(n)   &=& 
\sum_{i=0}^{N} (-1)^{i}  \si(\frac{n}{2^{i}}) \\
&=& \sum_{i=0}^{N} (-1)^{i} \left[ \frac{6}{\pi^2}\frac{n}{2^{i}} +O\left((\frac{n}{2^{i}})^{\theta}\right) \right] \\
&=& \frac{6n}{\pi^2}\sum_{i=0}^{N} (-\frac{1}{2})^{i} +\sum_{i=0}^{N} O\left((\frac{n}{2^{i}})^{\theta}\right) \\
&=& \frac{6n}{\pi^2} \frac{2}{3}\left(1-(- \frac{1}{2})^{N+1}\right)  +  n^{\theta}\ 
O\left (\sum_{i=0}^{N} \frac{1}{2^{i\theta}} \right) \\
&=& \frac{4n}{\pi^2} + O(1) +  n^{\theta}\ O(1-( \frac{1}{2^{\theta}})^{N+1}                 )   \\
&=& \frac{4n}{\pi^2} + O(n^{\theta})
\eeqas

Hence, by equation  (\ref{betti.shift})
$$ \sum_{k\ge 0} \be_k (\De_n) = \si^{odd} (n) - \si^{odd} (n/2) =
 \frac{4n}{\pi^2} + O(n^{\theta}) - 
\frac{2n}{\pi^2} - O((n/2)^{\theta}). $$
\eproof

\bthm\label{bettithm2} As $n\rarr \infty$ we have that
$$\sum_{k\ \mathrm{even}} \be_k (\De_n) \sim \frac{n}{\pi^2} \label{sumeven} 
\quad\mbox{ and }\quad
\sum_{k\ \mathrm{odd}} \be_k (\De_n) \sim \frac{n}{\pi^2} \label{sumodd} $$
\ethm
\bproof
Let $a(n)$ denote the first sum and $b(n)$ the second. Theorem \ref{bettithm1}
and the Prime Number Theorem show, respectively, that 
$$\frac{a(n)+b(n)}{n} \rarr \frac{2}{\pi^2} \quad\mbox{ and }\quad
\frac{a(n)-b(n)}{n} = \frac{M(n)}{n} \rarr 0$$
as $n \rarr\infty$.
Hence,
$$\frac{2a(n)}{n} = \frac{a(n)+b(n)}{n} + \frac{a(n)-b(n)}{n} \rarr \frac{2}{\pi^2}+0
$$
and similarly for $b(n)$.
\eproof

\bthm\label{bettithm3} For fixed $k$ and $n\rarr \infty$,
$$
\be_k(\De_n) \sim \frac{n}{2 \log n} \frac{(\log\log n)^k}{k!} \label{kvalue}
$$
\ethm

\bproof
We begin with a small auxiliary computation.
Using Theorem \ref{classics} we have that, as $n \rarr \infty$,

\beqas
\frac{\si_{k-i}(n/2^j)}{\si_{k}(n)} &\sim&
\frac{\frac{n/ 2^j}{\log (n/2^j)}}{\frac{n}{\log n}} \cdot 
\frac{\frac{(\log\log (n/2^j))^{k-i-1}}{(k-i-1)!}}{\frac{(\log\log n)^{k-1}}{(k-1)!}} \\
&=& \frac{1}{2^j} \cdot \frac{\log n}{\log (n/2^j)} \cdot 
\left( \frac{\log\log (n/2^j)}{\log\log n}   \right)^{k-i-1} \cdot \frac{(k-1)(k-2)\cdots (k-i)}{(\log\log n)^{i}} \\
&\rarr& \frac{1}{2^j} \cdot 1 \cdot 1 \cdot  \bca 1, \mbox{ if } i=0 \\
0,  \;\;\;\;\;\mbox{ if }  i > 0 
\eca
\eeqas
Hence, for fixed $i,j \ge 0$
\beq \lim_{n\rarr\infty} \frac{\si_{k-i}(n/2^j)}{\si_{k}(n)} = \bca 1/2^j, \mbox{ if } i=0 \\
0,  \;\;\;\;\;\mbox{ if }  i > 0 
\eca \label{lim}
\eeq
Equations (\ref{betti.shift}) and  (\ref{eq1}) imply
\beqas
\be_{k-1}(\De_n) &=&  \si^{odd}_{k} (n) - \si^{odd}_{k} (n/2) \\
&=& \sum_{j=0}^{k}(-1)^{j} \left[ \si_{k-j}(n/2^{j}) -  \si_{k-j}(n/2^{j+1})      \right] \\
\eeqas
Using  (\ref{lim}) we have that
\beqas
\frac{\be_{k-1}(\De_n)}{\si_k(n)} &=& 1-\frac{\si_k(n/2)}{\si_k(n)} +
\sum_{j=1}^{k}(-1)^{j} \left[ \frac{\si_{k-j}(n/2^{j})}{\si_k(n)} -  
\frac{\si_{k-j}(n/2^{j+1})}{\si_k(n)}      \right] \\
&\rarr& 1-\frac{1}{2} +0 = \frac{1}{2} 
\eeqas
Hence,
$$\be_{k-1}(\De_n) \sim  \frac{1}{2} \si_k(n) \sim  \frac{n}{2 \log n} \frac{(\log\log n)^{k-1}}{(k-1)!} 
$$
\eproof


\section{Some number-theoretic consequences}
The fact that the Betti numbers of the complex $\De_n$ are determined by 
counting certain odd squarefree integers has some purely number-theoretic implications.
Namely, 
knowing the 
number of such integers of weight $k+1$ in the interval $(0, n]$
one obtains  a lower bound for the
number of such integers of weight $k$ in the interval $(0, n/2]$.

We need to  recall the following   definitions.
Two number-theoretic functions  $\partial_k(n)$  and $\partial^k(n)$  are
defined in the following way. For \, $n,k\geq 1$ \, the integer $n$ can
in a unique way be expressed in the following form
$$ n=\binom {a_k}k +\binom {a_{k-1}}{k-1} +\cdots +\binom {a_i}i,$$

\noindent where \, $a_k>a_{k-1}>\cdots >a_i\geq i\geq 1.$ Then let:
$$ \partial_{k-1}(n)\defeq\binom {a_k}{k-1} +\binom {a_{k-1}}{k-2} +\cdots +\binom {a_i}{i-1},$$
and
$$ \partial^{k-1}(n)\defeq\binom {a_k-1}{k-1} +\binom {a_{k-1}-1}{k-2} +\cdots +\binom {a_i-1}{i-1}.$$

\noindent Also, we let $\partial_{k-1}(0)=\partial_{k-1}(0)=0.$

\bthm For all $k\ge 1$ we have that
\ben
\vspace{3mm}
\item  $\partial_{k}\left(\si_{k+1}^{odd}(n)\right) \le \si_{k}^{odd}(n/2)   $
\vspace{5mm}
\item  $\partial^{k}\left(\si_{2k+2}^{odd}(n)+\si_{2k+1}^{odd}(n)\right) \le \si_{2k}^{odd}(n/2)+\si_{2k-1}^{odd}(n/2)   $
\een
\ethm
\bproof
For an arbitrary finite simplicial complex with $f$-vector $(f_0, f_1, \ldots)$ and 
Betti numbers $\be_0, \be_1, \ldots$, let for $k\ge 0$
$$\chi_{k-1}\defeq\sum_{j\geq k} (-1)^{j-k}(f_j-\beta_j).$$ 
The following relations appear as 
Theorem 1.1 in \cite{BK} and Theorem 3.3 in \cite{BV}, respectively.
For all $k\geq 1$,

\beq\label{fbeta1}
 \partial_k(\chi_k+\beta_k)\leq \chi_{k-1}, \eeq


\beq\label{fbeta2}
 \partial^k(f_{2k+1}+\beta_{2k})\leq f_{2k-1}- \be_{2k-1}.\eeq

Let us now see what these relations mean for the particular complex $\De_n$.

\beqas
\chi_{k-1} &=& \sum_{j\ge k} (-1)^{j-k} \left(\si_{j+1}(n) -\be_{j} (\De_n)\right) 
= \sum_{j\ge k} (-1)^{j-k} \left( \si_{j+1}^{even}(n) + \si_{j+1}^{odd}(n/2)  \right) \\
&=& \sum_{j\ge k} (-1)^{j-k} \left( \si_{j}^{odd}(n/2) + \si_{j+1}^{odd}(n/2)  \right) 
= \si_{k}^{odd}(n/2).
\eeqas
Hence,
$$ \chi_{k}+\be_k= \si_{k+1}^{odd}(n/2)+\left( \si_{k+1}^{odd}(n) - \si_{k+1}^{odd}(n/2)\right)
= \si_{k+1}^{odd}(n),
$$
$$
\si_{2k+2}(n) +\be_{2k}(\De_n) = \si_{2k+2}(n) +\left( \si_{2k+1}^{odd}(n) - \si_{2k+1}^{odd} (n/2) \right)
= \si_{2k+2}^{odd}(n) +  \si_{2k+1}^{odd}(n),
$$
$$
\si_{2k}(n) -\be_{2k-1}(\De_n) = \si_{2k}(n) -\left( \si_{2k}^{odd}(n) - \si_{2k}^{odd} (n/2) \right)
= \si_{2k-1}^{odd}(n/2) +  \si_{2k}^{odd}(n/2).
$$
\vspace{.2mm}

\noindent
Inserting these evaluations into relations (\ref{fbeta1}) and (\ref{fbeta2}) one obtains the theorem.
\eproof

\section{The number-theoretic cell complex $\wt{\De}_n$}
\newcommand{\tdn}{\wt{\De}_n} 

The system of \sqf\ numbers less than or equal to $n$ and ordered by divisibility
presents itself immediately as a simplicial complex. The same is not true
for the larger system of {\em all} such numbers. However, a construction is known 
\cite{BV} which produces a CW complex 
of a similar nature. The CW structure
is no longer uniquely defined. However, if one demands that its closed $i$-dimensional 
cells are $(i-1)$-connected for all $i$ then the complex is 
uniquely determined up to homotopy type.

The construction from \cite{BV} of cellular realizations
of multicomplexes is reviewed in Section 6. We take it here for known.

The system of positive integers less than or equal to $n$, ordered by divisibility, is
isomorphic to a multicomplex of monomials. Namely, associate an indeterminate
$x_i$ with each prime number $p_i \le n$ and then extend this to a bijection
$$
p_{i_1}^{e_1}p_{i_2}^{e_2}\cdots  p_{i_k}^{e_k}, 
 \; \lrarr \;
x_{i_1}^{e_1}x_{i_2}^{e_2}\cdots  x_{i_k}^{e_k}, \quad i_1 < i_2 < \cdots < i_k$$
Therefore the construction of a cellular realization $\Ga$, reviewed in  Section 6, 
is applicable.

\begin{ddef}
Let $\tdn$ denote the CW complex  $\Ga(M)$, where $M$ is the multicomplex of
positive integers less than or equal to $n$, and the construction of
$\Ga$ is based on the choice of a well-connected CW string.
\end{ddef}

Here is a summary of the main features of $\tdn$, see Section 6 for further details.
\ben
\item The positive integers $k\le n$ are in bijection with the closed cells $c(k)$ of $\tdn$.
\item $\dim c(k)= \Om(k)-1$
\item The integer $k_1$ divides $k_2$ if and only if $c(k_1) \sbseq c(k_2)$.
\item 
The cell $c(k)$ has the following homotopy type:
$$c(k) \simeq \bca
S^{d-1}, \;\mbox{  if $k$ is a full square and  $\Om(k)=d$,} \\
\mbox{a point, \;otherwise}.
\eca
$$
\item In case we are using the $S^{\infty}$ string even more is true.
If $\Om (k)=d$, then the cell $c(k)$
has the following {\em homeomorphy} type (sphere or ball):
$$c(k) \cong \bca
S^{d-1}, \;\mbox{  if $k$ is a full square,} \\
B^{d-1}, \; \mbox{ otherwise}.
\eca
$$
\item The Euler characteristic is
$$\chi(\tdn) =  \# \{k \mid k\le n, \, \Om(k) \mbox{ is odd}\} -
\# \{k \mid k\le n, \, \Om(k) \mbox{ is even}\}
$$
\een

The Euler characteristic of $\tdn$ is equivalent to one of the classic functions of
number theory. The {\em summatory Liouville function} $L(n)$ is defined as follows.
$$ L(n)=\sum_{k\le n} (-1)^{\Om(k)}
$$
See \cite{BFM} for information about this function. So,
$$\chi(\tdn) = - L(n)
$$

Thus, the summatory Liouville function $L(n)$ plays for the cell complex
$\tdn$ the same role that the Mertens function $M(n)$ plays for
the simplicial complex $\De_n$. The two functions are related as follows.
\bprop\label{LvsM}
$$L(n)=\sum_{r=1}^{\lfloor\sqrt{n} \rfloor} M\left(\lfloor \frac{n}{ r^2} \rfloor \right)
\quad \mbox{ and } \quad
M(n)=\sum_{r=1}^{\lfloor\sqrt{n} \rfloor} \mu(r) \,L\left(\lfloor \frac{n}{ r^2} \rfloor \right)
$$
\eprop
\bproof
Let $[n]\defeq \{1,2, \ldots ,n\}$ and $[n]_{\mathrm{sqf}} \defeq \{k\in [n] \mid k \mbox{ is \sqf}\} $.
Every positive  integer $k$ can be uniquely factored as a product of a full square
$r^2$ and a \sqf\ integer $s$: $k=r^2 s$. This implies the following
set-theoretic disjoint union decomposition
$$[n]= \biguplus_{r=1}^{\lfloor\sqrt{n} \rfloor} r^2 \cdot  \left[ \lfloor \frac{n}{ r^2} \rfloor 
\right]_{\mathrm{sqf}}
$$
Since $\Om(k)=\Om(r^2 s) \equiv \Om(s)\, (\mathrm{mod} 2)$, the first formula follows.
Then the second one is obtained by M\"obius inversion of the form
presented in \cite[p. 307]{HW}.
\eproof

The rates of growth  of the two functions $L(n)$ and $M(n)$ are essentially
identical, as the following proposition shows.

\bprop\label{LMequiv}
For every $\eps> 0$, as $n\rarr\infty$,
$$M(n) =o(n^{1/2+\eps}) \mbox{\; if and only if \;} L(n) =o(n^{1/2+\eps})
$$
\eprop
\bproof In what follows all numbers that are not integers are to be rounded down to the 
closest smaller integer. So, $\sqrt{n}$ is to be read $\lfloor\sqrt{n}\rfloor$, etc.

Let $\eps>0$ and assume that $M(n) =o(n^{1/2+\eps})$ as $n \rarr \infty$.
Choose any $\de>0$, and let $N$ be such that
$\sum_{r=N+1}^{\infty} \frac{1}{r^{1+2\eps}} <\de$ and $1/N \le \de$. Then, if
$n$ is large enough that $\sqrt{n}> N$ and
$$ \frac{|M(\frac{n}{r^2})|}{(\frac{n}{r^2})^{1/2+\eps}} \le 1/ N^2, \mbox{ for all } r \in \{1,2, \cdots, N\},
$$
we have that
\beqas \frac{[L(n)|}{n^{1/2+\eps}} &\le & 
\sum_{r=1}^{\sqrt{n}}
\frac{1}{r^{1+2\eps}} \; \frac{|M(\frac{n}{r^2})|}{(\frac{n}{r^2})^{1/2+\eps}} \\
&\le&
\sum_{r=1}^{N}
\frac{1}{r^{1+2\eps}} \; \frac{|M(\frac{n}{r^2})|}{(\frac{n}{r^2})^{1/2+\eps}} + \sum_{r=N+1}^{\infty}
\frac{1}{r^{1+2\eps}} \\
&\le& N\cdot \frac{1}{N^2} +\de \\
&\le&  2\de.
\eeqas

The same computation, exchanging the roles of $M(n)$ and $L(n)$, gives by
Proposition \ref{LvsM} the opposite implication.
\eproof

By way of the $\eps=1/2$ case of Proposition
\ref{LMequiv}, the Prime Number Theorem  implies that
\beq\label{altPNT}
L(n)= o(n), \; \mbox {as } n\rarr\infty.
\eeq
Also, one can conclude from Proposition \ref{LMequiv} that
the Riemann Hypothesis is equivalent to 
\beq\label{altRH}
L(n)= O(n^{1/2+\eps}), \; \mbox {as } n\rarr\infty, \mbox { for every }
 \eps > 0.
\eeq

The preceding shows that the Euler characteristic of the cell complex 
$\tdn$ has the same number-theoretic relevance as that of the
simplicial complex $\De_n$. This motivates seeking information about
the Betti numbers of $\tdn$, as we did in Section 3 for $\De_n$.

\bthm\label{homotopy2}
We have the following homotopy equivalence
$$\tdn \simeq \bigvee_{r^2\le n} \mathrm{ susp}^{2\Om(r)}  (\De_{\lfloor n/r^2 \rfloor})
$$
\ethm
\bproof
This is a direct consequence of Theorem \ref{homotopy}.
\eproof

\bcor\label{homsum}
$\tdn$ has the homotopy type of a wedge of spheres, and
$$\be_k(\tdn)= \sum_{r=1}^{\lfloor\sqrt{n}\rfloor} \be_{k-2\Om(r)}  (\De_{\lfloor n/r^2 \rfloor})
$$
\ecor
\bproof
Combine with Theorem \ref{homotopy2} the information coming from
Theorems \ref{betti1} and \ref{hom}.
\eproof

\bthm\label{bettithm5} As $n\rarr \infty$ we have that
$$ \sum_{k\ge 0} \be_k (\tdn) = \frac{n}{3} +O(n^{\theta}),\;\;  \mbox{ for all } \theta > \frac{22}{27}$$
\ethm
\bproof Let $\psi> \frac{17}{54}$. Using Theorem \ref{bettithm1} and Corollary \ref{homsum} we get:
\beqas  \sum_{k\ge 0} \be_k (\tdn)  &=& \sum_{r=1}^{\lfloor\sqrt{n}\rfloor} \sum_{k} \be_k (\De_{\lfloor n/r^2 \rfloor}) \\
&=& \sum_{r=1}^{\lfloor\sqrt{n}\rfloor} \left(\frac{2(n/r^2)}{\pi^2} + O((\frac{n}{r^2})^{\psi}) \right) \\
&=& \frac{2n}{\pi^2}\sum_{r=1}^{\lfloor\sqrt{n}\rfloor} \frac{1}{r^2}+ \sqrt{n}\cdot  O(n^{\psi}) \\
&=& \frac{2n}{\pi^2} \frac{\pi^2}{6} + O(n \sum_{r=\lfloor\sqrt{n}\rfloor}^{\infty} \frac{1}{r^2})+ O(n^{1/2+\psi}) \\
&=& \frac{n}{3} + O(n^{1/2})+ O(n^{1/2+\psi}) \\
&=& \frac{n}{3} + O(n^{\theta}), \mbox{ where } \theta =\frac{1}{2}+\psi > 
\frac{1}{2}+ \frac{17}{54}= \frac{22}{27}. 
\eeqas
\eproof

\bthm\label{bettithm6} As $n\rarr \infty$ we have that
$$\sum_{k\ \mathrm{even}} \be_k (\tdn) \sim \frac{n}{6}  
\quad\mbox{ and }\quad
\sum_{k\ \mathrm{odd}} \be_k (\tdn) \sim \frac{n}{6}  $$
\ethm
\bproof 
Let $a(n)$ denote the first sum and $b(n)$ the second. Theorem \ref{bettithm5}
and equation (\ref{altPNT}) show, respectively, that 
$$\frac{a(n)+b(n)}{n} \rarr \frac{1}{3} \quad\mbox{ and }\quad
\frac{a(n)-b(n)}{n} = \frac{L(n)}{n} \rarr 0$$
as $n \rarr\infty$.
Hence, the proof can be completed as that of Theorem \ref{bettithm2}.

\eproof

\section{Cellular realization of multicomplexes}

In this section we review the construction of CW complexes from \cite{BV}.
We expand on some details and and present
a new result concerning their homotopy type.

By a {\em CW string} we shall mean an infinite sequence of embeddings
$$c^0 \hookrightarrow c^1 \hookrightarrow  \cdots \hookrightarrow c^j \hookrightarrow \cdots
$$
such that
\ben 
\item $\{c^j\}_{j\ge 0}$ are the closed cells of a CW decomposition of the colimit $\cup c^j$.
\item $\dim c^j =j$, for all $j\ge 0$.
\een
The CW string is said to be {\em well-connected} if 
\ben 
\item[(3)] $c^j$ is $(j-1)$-connected, for all $j\ge 0$.
\een

\blem
A CW string is well-connected if and only if
$c^j $ is contractible for all even $j$ and has the homotopy type of the $j$-sphere for
all odd $j$.
\elem
\bproof
It is known that a $(j-1)$-connected and $j$-dimensional space has the
homotopy type of a wedge of $j$-dimensional spheres. The number of spheres 
is given by the Euler characteristic. Since the string has exactly
one cell in each dimension, the reduced Euler characteristic of $c^j$ is $0$ 
for all even $j$ and $1$ for all odd $j$.
\eproof

Here are three constructions of CW strings. The first one was suggested to us by
M. Falk (private communication). The other two appear in \cite{BV}.
\vspace{5mm}

1. The $S^{\infty}$ string.
This is a CW decomposition of the infinite-dimensional
 sphere $S^{\infty}=\{(x_1, x_2, \ldots)\mid x_n =0 \mbox{ for sufficiently large $n$, and }
x_1^2+ x_2^2+ \cdots =1\}$. The cells are
\beqas
c^j &=& \{(x_1, x_2, \ldots) \sbseq S^{\infty} \mid
x_k =0 \mbox{ for all }k > j+1\}, \mbox{ if $j$ is odd,} \\
c^j &=& \{(x_1, x_2, \ldots) \sbseq S^{\infty} \mid
x_{j+1} \ge 0, \; x_k =0 \mbox{ for all }k > j+1\}, \mbox{ if $j$ is even.}
\eeqas

It is clear, when $j$ is odd, how to attach a $(j+1)$-ball onto the sphere $c^j$  in order to obtain 
the ball $c^{j+1}$. When $j$ is even the attachment map has the following
description. Note that $c^j$ is here the upper hemisphere of $S^{j}$, the $j$-dimensional
sphere. Consider the $(j+1)$-ball
$$ B^{j+1} = \{(x_1, x_2, \ldots) \sbseq S^{\infty} \mid
x_{j+2} \ge 0, \; x_k =0 \mbox{ for all }k > j+2\},
$$
i.e., the upper hemisphere of $S^{j+1}$. By switching to polar coordinates for the
last two coordinates, the points of $S^{j}$ can be coordinatized
$$ S^{j} = \{(x_1, x_2, \ldots, x_{j-1}, r\cos \theta, r \sin \theta, 0, 0, \ldots)\} \sbseq S^{\infty}$$

In terms of these coordinates the attachment map $\al: \mbox{bd}(B^{j+1}) = S^j \rarr c^j$
sends  $(x_1, x_2, \ldots, x_{j-1}, r\cos \theta, r \sin \theta)$ to 
$(x_1, x_2, \ldots, x_{j-1}, r\cos 2 \theta, r \sin 2 \theta)$.

Thus we have a CW string which is well-connected in a strong sense, namely
its cells are {\em homeomorphic} to balls and spheres of the appropriate dimensions.

\vspace{5mm}
 
 2. The $D^{\infty}$ string.
The {\em $j$-dimensional dunce hat} $D^j$ is obtained from the $j$-dimensional
simplex $(a_0, a_1, \ldots, a_j)$ by identifying all of its  $(j-1)$-dimensional
faces, with their vertices in the  induced order. Note that $D^2$ is the usual dunce hat. It is
shown in \cite{AMS} that $D^j$ is contractible, but not collapsible, for all
even $j>0$, and that $D^j$ is a homotopy sphere for all odd $j$.
Thus, with the obvious attaching maps we have a well-connected CW string 
$$D^0 \hookrightarrow D^1 \hookrightarrow  \cdots \hookrightarrow D^j \hookrightarrow \cdots
$$
giving a CW decomposition of its colimit $D^{\infty}$, the infinite-dimensional
dunce hat.

\vspace{5mm}

3. The $\R P^{\infty}$ string.
The standard CW decomposition of infinite-dimensional real projective space, with one
cell in each dimension, provides a well-known example of a CW string
$$\R P^{0} \hookrightarrow \R P^{1} \hookrightarrow  \cdots \hookrightarrow 
\R P^{j} \hookrightarrow \cdots
$$
This string is not well-connected. However, on the level of rational homology
$H(\cdot, \Q)$ this string behaves in a way that for algebraic purposes parallels
that of well-connected strings, see Remark \ref{realproj}.

\vspace{5mm}

Let $M$ be a multicomplex.
By this we mean a finite collection of monomials in indeterminates
$x_1, x_2, \ldots$ closed under divisibility. A CW complex $\Ga(M)$ is constructed as follows.
It depends on a choice of CW string $\{c^j\}_{j\ge 0}$, which once chosen remains fixed and will not be 
included in the notation. 

For each indeterminate $x_i$ take a copy $\{c^j_i\}_{j\ge 0}$
of the string. Then, to each monomial $m=x_{i_1}^{e_1} x_{i_2}^{e_2} \cdots x_{i_k}^{e_k}\in M$ associate the space
$$c(m)\defeq c_{i_1}^{e_1 -1} \ast c_{i_2}^{e_2-1} \ast \cdots\ast  c_{i_k}^{e_k-1}
$$
where ``$\ast$'' denotes the join of topological spaces. Note that if $m$ is squarefree
then $c(m)$ is a $(k-1)$-dimensional simplex, since $c^0$ is a single  vertex
(for every choice of CW string).

We say that $m$ 
is a {\em full square} if $m=r^2$ for some monomial $r$.

\blem \cite[Prop. 2.2]{BV}\label{cells}
The space $c(m)$ has the following homotopy type:
$$c(m) \simeq \bca
S^{d-1}, \;\mbox{  if $m$ is a full square of degree $d$,} \\
\mbox{a point, \;otherwise}.
\eca
$$
\elem
We remark that in case we are using the $S^{\infty}$ string even more is true: $c(m)$ is
{\em homeomorphic} to the $(d-1)$-sphere or the $(d-1)$-ball, respectively.

The CW complex $\Ga(M)$ associated with $M$ is constructed as follows. With each
monomial $m\in M$ associate the space $c(m)$, and then glue these together to form
$\Ga(M)$. The attaching maps are everywhere the ones coming from how a $j$-cell 
is attached to $c^{j-1}_i$ to obtain $c^j_i$, and the joins of these maps.

\bprop \cite[Sect. 2]{BV}
The space $\Ga(M)$ is a CW complex with closed cells $c(M)$, $m\in M$.
A cell $c(m)$ is contained in another cell $c(m')$ if and only if $m$ divides $m'$.
\eprop

We think of $\Ga(M)$ as the geometric realization of the multicomplex $M$.
Note that if all monomials in  $M$  are squarefree then the construction reduces
to the usual geometric realization of a simplicial complex.

Every monomial can be uniquely factorized as a product of a full square monomial
and a \sqf\ monomial: $m=r^2 s$. For each full square monomial $r^2$ in $M$ we define
$$M_{r^2} = \{s\in M \mid r^2 s\in M, s \; \mbox{\sqf} \}
$$
Notice that  each $M_{r^2}$ is a simplicial complex, possibly empty. 

Let $|r|$ denote the degree of a monomial $r$.

\bthm \label{homotopy}
Suppose that a well-connected CW string has been used for the construction of  $\Ga(M)$. Then
$$\Ga(M) \simeq \bigvee_{r^2\in M} \mathrm{ susp}^{2|r|}\,  \Ga(M_{r^2})
$$
\ethm
\bproof
We are going to use the theory of homotopy colimits of diagrams of spaces.
The tools from this theory that we need are summarised in an
accessible way in \cite{WZZ}. We
refer the reader to that paper for all explanations of terminology and 
basic facts used in the sequel. We frequently use the fact that the join operation of
topological spaces commutes with the needed operations, up to homotopy.

The multicomplex $M$, ordered by divisibility is a poset. The functor
$$\DD: m\mapsto c(m) \quad  \mathrm{ and }\quad \DD: m | m' \mapsto \left(c(m) \hookrightarrow c(m')\right)$$
gives us a diagram of spaces $\DD: M \rarr Top$. By the ``Projection Lemma'' \cite[Prop. 3.1]{WZZ}
we have that
\beq\label{diag1}
\Ga(M) \simeq \mathrm{hocolim}\, \DD
\eeq

Let $\EE$ be a diagram over $M$ with constant maps. Since $\dim\left( \bigcup_{m'<m} c(m')\right)
= \dim c(m)-1$ and $c(m)$ is $(\dim c(m)-1)$-connected, it follows that any two maps
from $\bigcup_{m'<m} c(m')$ to $c(m)$ are homotopic. In particular, the inclusion map
is homotopic to any constant map. Therefore, using the ``Homotopy Lemma'' \cite[Prop. 3.7]{WZZ}
and the homotopy extension property, we get
\beq\label{diag2}
\mathrm{hocolim\, \DD} \simeq \mathrm{hocolim\, \EE}
\eeq

Finally, the ``Wedge Lemma'' \cite[Lemma 4.9]{WZZ} implies that
\beq\label{diag3}
\mathrm{hocolim\, \EE} \simeq  \bigvee_{m\in M}  \left( c(m) \ast \De(M^{>m})\right),
\eeq
where $\De(M^{> m})$ denotes the order complex of the poset of monomials in $M$
that are strictly above $m$ in the partial order, i.e., that are  divisible by $m$. 

Many of the spaces appearing in the wedge are contractible and can
therefore  be removed without affecting the homotopy type. Namely,
we have that
$$ c(m) \ast \De(M^{>m}) \; \simeq \; \bca
\mathrm{susp}^{2|r|}(\De(M^{>m})), \; \mathrm{if}\, m=r^2\in M \\
\mbox{a point, \quad\qquad otherwise}
\eca
$$
This follows from Lemma \ref{cells} and the fact that taking join with the
$(d-1)$-dimensional sphere is, up to homeomorphism, equivalent to
taking $d$-fold suspension.
Combining equations (\ref{diag1}), (\ref{diag2}) and (\ref{diag3}) we obtain the homotopy
equivalence
\beq\label{diag4}
\Ga(M) \simeq  \bigvee_{r^2 \in M}  \mathrm{susp}^{2|r|}\De(M^{>r^2})
\eeq
Thus, the following homotopy equivalences  complete the proof:
$$\De(M^{>r^2}) \simeq \De(M_{r^2}) \simeq  \Ga(M_{r^2})
$$
The first equivalence  is implied by \cite[Prop. 4.1]{BV}. The second 
follows from the fact that the homeomorphism type of a
simplicial complex is invariant under barycentric subdivision.

\eproof

\bcor\cite[p.55]{BV}\label{hom}
$$ \wt{H}_k (\Ga(M),\Z) \cong\bigoplus_{r^2\in M} \wt{H}_{k-2|r|} (\Ga(M_{r^2}), \Z)
$$
\ecor

\brem\label{realproj} \rm{
The  $\R P^{\infty}$ string is not well-connected, since its cells $c^j$ are real projective
spaces with plenty of mod $2$ torsion. However, the string is ``rationally well-connected'',
 since
$$\wt{H}_k (c^j ,\Q) \cong \bca \Q, \mbox{if $j$ is odd and $k=j$,} \\
0, \mbox{in all other cases}.
\eca $$
Thus, in terms of rational  homology 
the $\R P^{\infty}$ string behaves  like the well-connected strings.
By replacing the diagrammatic  tools used in the proof of Theorem \ref{homotopy}
by their homological counterparts one can derive a splitting formula for
$ \wt{H}_k (\Ga(M),\Q)$ analogous to Corollary \ref{hom}. The rational Betti
numbers produced for the number-theoretic cell complexes using joins of 
real projective spaces as cells will therefore be the same as those
computed in Sections 3 and 5.
}
\erem


\begin{thebibliography}{}
\bibitem[AMS]{AMS} R.~N.~Andersen, M.~M.~Marjanovi\'c, R.~M.~Shori,
Symmetric products and higher-dimensional dunce hats,
{\em Topology Proc.} {\bf 18} (1993), 7--17.

\bibitem[BK]{BK} A.~Bj\"orner, G.~Kalai, An extended Euler-Poincar\'e theorem,
{\em Acta Math.} {\bf 161} (1988), 279--303.

\bibitem[BV]{BV} A.~Bj\"orner, S.~Vre\'cica, On $f$-vectors and Betti numbers of multicomplexes,
{\em Combinatorica} {\bf 17} (1997), 53--65.

\bibitem[BFM]{BFM} P.~Borwein, R.~Ferguson, M.~J.~Mossinghoff, 
Sign changes in sums of the Liouville function, {\em Math. Comp.}
{\bf 77} (2008), 1681--1694.

\bibitem[HW]{HW} G.~H.~Hardy, E.~M.~Wright, {\em An Introduction to the Theory
of Numbers}, Sixth Edition, Oxford Univ. Press, Oxford, UK, 2008.

\bibitem[Iv]{Iv} A.~Ivi\'c, {\em The Riemann zeta-function}, Wiley-Interscience, New York, 1985.

\bibitem[WZZ]{WZZ} V.~Welker, G.~M.~Ziegler, R.~T.~Zivaljevi\'c, Homotopy colimits --- comparison 
lemmas for combinatorial applications,
{\em J. reine angew. Math.} {\bf 509} (1999), 117--149.

\end{thebibliography}
\end{document}